\documentclass[USenglish,cleveref]{lipics-v2021}

\title{A Sublinear Bound on the Page Number of Upward Planar Graphs\thanks{A preliminary version of this paper is published in the Proceedings of the 2022 Annual \textsc{acm-siam} Symposium on Discrete Algorithms (\textsc{soda})~\cite{JMU_UpwardPageNumber}}}
\hypersetup{
  pdftitle={A Sublinear Bound on the Page Number of Upward Planar Graphs},
}

\author{Paul Jungeblut}%
  {Karlsruhe Institute of Technology, Karlsruhe, Germany}%
  {paul.jungeblut@kit.edu}
  {https://orcid.org/0000-0001-8241-2102}{}

\author{Laura Merker}%
  {Karlsruhe Institute of Technology, Karlsruhe, Germany}%
  {laura.merker2@kit.edu}%
  {https://orcid.org/0000-0003-1961-4531}{}%

\author{Torsten Ueckerdt}%
  {Karlsruhe Institute of Technology, Karlsruhe, Germany}%
  {torsten.ueckerdt@kit.edu}%
  {}{}%

\authorrunning{P. Jungeblut, L. Merker, and T. Ueckerdt}


\nolinenumbers
\hideLIPIcs

\graphicspath{{fig/}}
\makeatletter
\def\input@path{{content/}}
\makeatother

\usepackage{mathtools}
\usepackage{xspace}
\usepackage{csquotes}
\usepackage{microtype}
\usepackage{cite}

\DeclareMathOperator{\pn}{pn}
\DeclareMathOperator{\tn}{tn}

\newcommand{\calP}{\ensuremath{\mathcal{P}}\xspace}
\newcommand{\calQ}{\ensuremath{\mathcal{Q}}\xspace}
\renewcommand{\leq}{\leqslant}
\renewcommand{\geq}{\geqslant}
\renewcommand{\preceq}{\preccurlyeq}
\renewcommand{\succeq}{\succcurlyeq}
\DeclareMathOperator{\Grid}{Grid}
\newcommand{\gedge}[2]{\ensuremath{\bigl(#1,#2\bigr)}\xspace}

\begin{document}

\maketitle

\begin{abstract}
  The page number of a directed acyclic graph~$G$ is the minimum~$k$ for which there is a topological ordering of~$G$ and a~$k$\nobreakdash-coloring of the edges such that no two edges of the same color cross, i.e., have alternating endpoints along the topological ordering.
  We address the long-standing open problem asking for the largest page number among all upward planar graphs.
  We improve the best known lower bound to~$5$ and present the first asymptotic improvement over the trivial~$\mathcal{O}(n)$ upper bound, where~$n$ denotes the number of vertices in~$G$.
  Specifically, we first prove that the page number of every upward planar graph is bounded in terms of its width, as well as its height.
  We then combine both approaches to show that every $n$-vertex upward planar graph has page number $\mathcal{O}(n^{2/3} \log^{2/3}(n))$.
\end{abstract}

\section{Introduction}
\label{sec:introduction}

\begin{figure}[t]
  \centering
  \includegraphics{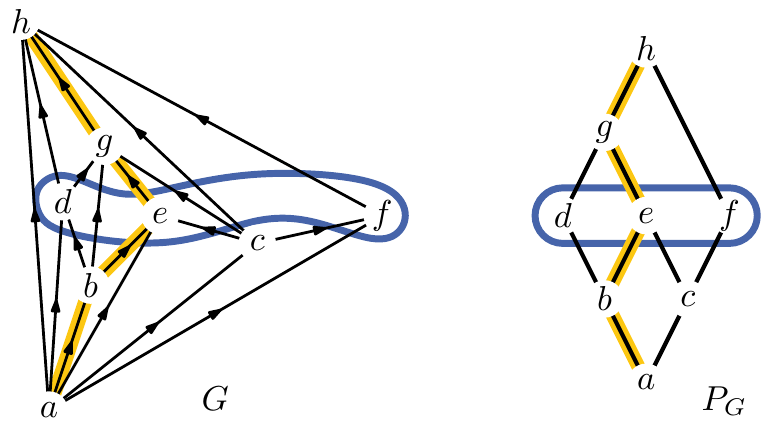}
  \hspace{3em}
  \includegraphics{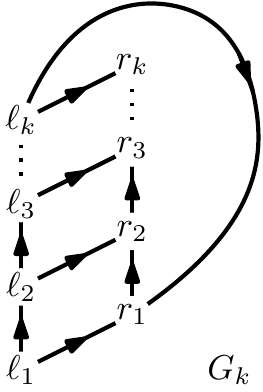}
  \caption{
    Left: An upward planar $st$-graph $G$ of height~$5$ and width~$3$.
    Middle: The reachability poset $P_G$ of $G$.
    Right: A planar directed acyclic graph $G_k$ with $\pn(G_k) = k$.}
 \label{fig:examples}
\end{figure}

In an \emph{upward planar drawing} of a directed acyclic graph $G = (V,E)$, every vertex $v \in V$ is a point in the Euclidean plane, and every edge $(u,v) \in E$ is a strictly $y$-monotone curve\footnote{Equivalently, straightline segments may be used~\cite{BT88}.} with lower endpoint $u$ and upper endpoint $v$ that is disjoint from other points and curves, except in its endpoints.
A directed acyclic graph admitting such a drawing is called \emph{upward planar}.
In other words, a directed graph is upward planar if it allows a planar drawing with all edges ``going strictly upwards''.
In \cref{fig:examples} we have an upward planar graph $G$ on the left, while the planar directed acyclic graph $G_k$ on the right is not upward planar.

In a \emph{book embedding} of a directed acyclic graph~$G = (V,E)$, the vertex set~$V$ is endowed with a topological ordering~$ < $, called the \emph{spine ordering}, and the edge set~$E$ is partitioned into so-called \emph{pages}
with the property that no page contains two edges $(u_1,v_1), (u_2,v_2)$ that cross with respect to $ < $, i.e., $u_1 < u_2 < v_1 < v_2$ or $u_2 < u_1 < v_2 < v_1$.
Then the \emph{page number}~$\pn(G)$ of a directed acyclic graph~$G$ is the minimum~$k$ for which it admits a book embedding with~$k$ pages.
In other words, $\pn(G) \leq k$ if the vertices can be ordered along the (horizontal) spine with all ``edges going right'' and there exists a~$k$\nobreakdash-edge coloring so that any two edges with alternating endpoints along the spine have distinct colors.

\begin{figure}
  \centering
  \includegraphics{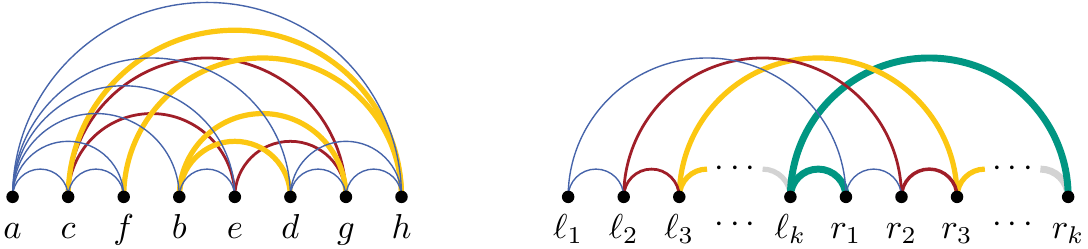}
  \caption{
    Book embeddings of the graphs in \cref{fig:examples}.
  }
  \label{fig:book-embeddings}
\end{figure}

In \cref{fig:book-embeddings} we have book embeddings of the directed acyclic graphs $ G $ and $ G_k $ in \cref{fig:examples} with three pages (left) and~$k$ pages (right), respectively.
This shows that $\pn(G) \leq 3$ and $\pn(G_k) \leq k$.
In fact, observe that $G_k$ admits only one topological ordering $<$, as there is a directed Hamiltonian path $\ell_1,\ldots,\ell_k,r_1,\ldots,r_k$ in $G_k$.
As the edges $(\ell_1,r_1),\ldots,(\ell_k,r_k)$ are pairwise crossing w.r.t.~$<$, it follows that $\pn(G_k)=k$.
(Recall that $G_k$ from the right of \cref{fig:examples} is not upward planar.)
It is easy to see (as observed for example in~\cite{BDDDMP19}) that for any directed graph $G$ we have $\pn(G)\leq 2$ if and only if $G$ is a spanning subgraph of an upward planar graph with a directed Hamiltonian path.
It thus follows that $\pn(G)=3$ for the graph $G$ in the left of \cref{fig:examples}.

The page number of \emph{undirected} graphs (where the spine ordering may be any vertex ordering) was introduced by Bernhart and Kainen in 1979~\cite{BK79}, building upon the suggested notion of Ollmann~\cite{Oll73}.
Their conjecture that the page number of planar graphs is unbounded was quickly disproven~\cite{BS84,Hea84}, with Yannakakis~\cite{Yan89} giving the best upper bound of $4$, which was just very recently shown to be best-possible~\cite{Yan20,BKKPRU20}.

Book embeddings of \emph{directed} graphs were first considered by Nowakowski and Parker~\cite{NP89} in 1989.
They introduce the page number of a poset $P$ by considering its cover graph $G(P)$ and restricting the spine ordering to be a topological ordering of $G(P)$, or equivalently, a linear extension of $P$.
They then ask whether posets with a planar order diagram have bounded page number --- equivalently, whether upward planar and transitively reduced graphs have bounded page number.
Despite significant effort on posets~\cite{MJZ15,AR96,ARK01,HP99,HPT99,HP97,Hun93,Sys89} and upward planar graphs~\cite{GLMSW15,DDGW06,MS09,FFR13,ADD13,BDDDMP19,NP21,BDMN21,BDFGMR22}, this question is still open.
In fact, the question was not only asked for upward planar graphs but also for several subclasses thereof.
Apart from planar posets, these include upward outerplanar graphs~\cite{NP21,BDMN21}, upward planar 2-trees~\cite{NP21,GDLW06}, and planar lattices~\cite{NP89}.
For all these graph classes, the asymptotically best known upper bound is linear in $n$, the number of vertices, which can be obtained by simply putting each edge on a separate page.


In this paper, we provide the first upper bound for any upward planar graph that is sublinear in the number of vertices.
Specifically, we prove that~$n$-vertex upward planar graphs have page number~$ \mathcal{O}(n^{2/3} \log^{2/3}(n))$.
We do so by bounding the page number of any upward planar graph first in terms of its width, then in terms of its height, and finally combining both approaches to achieve the desired bound in terms of its number of vertices.

\subparagraph{Related Work.}
Nowakowski and Parker~\cite{NP89} (and independently Heath et al.~\cite{HPT99}) show that directed forests have page number $1$.
Alzohairi and Rival~\cite{AR96} (see also~\cite{DDGW06}) show that series-parallel digraphs have page number $2$, which was later generalized to $N$-free upward planar graphs by Mchedlidze and Symvonis~\cite{MS09}.

The best known upper bounds for the page number of upward planar graphs are due to
Frati et al.~\cite{FFR13}, who prove that
every $n$-vertex upward planar triangulation with $o(n/\log n)$ diameter has $o(n)$ page number and that
every $n$-vertex upward planar triangulation has page number $o(n)$ if that is true for those with maximum degree $\mathcal{O}(\sqrt{n})$.
They further prove that any $n$-vertex upward planar triangulation has page number at most $\min\{\mathcal{O}(k \log n), \mathcal{O}(2^k)\}$, where $k$ is the maximum page number among its $4$-connected subgraphs.
This bound is improved by a result by Davies~\cite{Dav22} to $ \mathcal{O}(k \log k) $.
According to the authors of~\cite{FFR13} \enquote{Determining whether every $n$-vertex upward planar DAG has $o(n)$ page number [\ldots] remains among the most important problems in the theory of linear graph layouts.}

For lower bounds, Nowakowski and Parker~\cite{NP89} present a planar poset with page number~$3$, while Hung~\cite{Hun93} presents a planar poset with page number~$4$.
For general upward planar graphs one can also easily derive the same lower bound of $4$ from one of the undirected planar graphs of page number $4$~\cite{Yan20,BKKPRU20}.
As nothing better is known here, we also present in this paper upward planar graphs with page number at least $5$.

\subparagraph{Preliminaries.}
We denote the directed reachability of a vertex~$v$ from another vertex~$u$ in a directed acyclic graph~$G$ by $u \prec_G v$ (omitting the index if it is clear from the context), and write $u \preceq v$ if $u \prec v$ or $u=v$.
This way we obtain the \emph{reachability poset} $P_G = (V,\prec)$ of $G$ as the vertices of $G$ partially ordered by their directed reachability.
Transferring these notions from posets to directed acyclic graphs, we say that $ u $ and $ v $ are \emph{comparable} if $ u \preceq v $ or $ v \preceq u $; otherwise $u$ and $v$ are \emph{incomparable}.
Consequently, the \emph{height}~$h(G)$ and \emph{width}~$w(G)$ of a directed acyclic graph~$G$ is the largest number of pairwise comparable, respectively incomparable, vertices in $G$.
Equivalently, $h(G)$ is the number of vertices in a longest directed path in~$G$, while~$w(G)$ is the largest number of vertices in~$G$ with no directed reachabilities among them.
See the left and middle of \cref{fig:examples} for some example.
Let us also define for a subset $X$ of vertices of $G$ its height $h(X)$ and width $w(X)$ as the maximum number of vertices in $X$ that are pairwise comparable, respectively incomparable, in $G$.

An upward planar graph $G = (V,E)$ is an \emph{$st$-graph} if there is a (unique) vertex~$s$ with $s \preceq v$ for all $v \in V$ and a (unique) vertex $t$ with $v \preceq t$ for all $v \in V$.
An $st$-path in $G$ is a directed path from $s$ to $t$ in $G$.
In particular, the height of an $st$-graph is the length of a longest $st$-path.
It is known~\cite{Kel87} that every upward planar graph~$G$ (on at least three vertices) is a spanning subgraph of some $st$-graph~$\overline{G}$ whose faces are all bounded by triangles.
Note that this augmentation might not be unique.
As $\pn(G) \leq \pn(\overline{G})$ whenever $G \subseteq \overline{G}$, we may restrict ourselves to $st$-graphs when proving upper bounds on the page number of upward planar graphs in terms of their number of vertices.
(Note however that this is not true when working in terms of height.)
Let us also remark that if~$G$ is an $st$-graph, its reachability poset~$P_G$ is called a \emph{planar lattice} in order theory~\cite{BFR72}.

A notion closely related to the page number is the \emph{twist number}~$\tn(G)$, which is defined as the smallest~$k$ for which there exists a topological ordering~$<$ of~$G$ with no~$(k+1)$\nobreakdash-twist, i.e., no~$k+1$ edges that are pairwise crossing w.r.t.\ $<$.
Clearly, $\tn(G) \leq \pn(G)$, as the $k$ edges of a $k$-twist must be assigned to pairwise distinct pages.
Indeed, having already decided on a spine ordering with no $(k+1)$-twist, assigning the edges to pages is equivalent to coloring the vertices of a corresponding circle graph $H$ with no $(k+1)$-clique.
Since circle graphs are $\chi$-bounded~\cite{Gya85}, one can actually bound the number of pages in terms of the largest twist size.
The currently best result due to Davies~\cite{Dav22} states that
$ \chi(H) \leq 2 \omega(H) \log_2 \omega(H) + 2 \omega(H) \log_2 \log_2 \omega(H) + 10 \omega(H) \leq 14 \omega(H) \log_2 \omega(H) $
for every circle graph~$H$ (where~$\omega(H)$ is the clique number of~$H$), which gives the following.%

\begin{observation}\label{obs:twist-to-page}
 For every directed acyclic graph $G$ we have $\pn(G) \leq 14 \tn(G) \log_2 \tn(G) $.
 Moreover, the bound can be certified by a book embedding using any vertex ordering with maximum twist size $ \tn(G) $.
\end{observation}

%

In fact, we shall often times bound the twist number of the considered upward planar graph $G$ and then conclude for its page number via \cref{obs:twist-to-page}.

  All graphs considered in this paper are directed and in most figures we omit the arrows indicating an edge's direction.
  If not explicitly drawn otherwise, all edges are oriented upwards.

\subparagraph{Our Results.}
First, we bound the page number of upward planar graphs $G$ in terms of their width.

\begin{theorem}\label{main:width}
 Every upward planar graph $G$ of width $w$ has $\pn(G) \leq 14 w$.
\end{theorem}

Then, we bound the page number of $st$-graphs in terms of their height.
In fact, we show that $\tn(G) \leq 4h(G)$, improving on the $\tn(G) \leq \mathcal{O}( h(G)\log(n)) $ bound for every~$n$\nobreakdash-vertex $st$-graph~$G$ due to Frati et al.~\cite{FFR13}.\footnote{Frati et al.~\cite{FFR13} refer to $h(G)$ as the diameter of $G$.}
Together with \cref{obs:twist-to-page} this gives the following.

\begin{theorem}\label{main:height}
 Every upward planar graph $G$ of height $h$ has $\pn(G) \leq 56 h (\log h + 2) $.
\end{theorem}

Combining our approaches for bounded width and bounded height, we give the first sublinear upper bound on $\pn(G)$ in terms of the number of vertices in $G$.

\begin{theorem}\label{main:vertices}
 Every upward planar graph $G$ on $n$ vertices has $ \pn(G) \leq \mathcal{O}(n^{2/3} \log^{2/3}(n)) $.
\end{theorem}

Finally, we improve the best known lower bound on the maximum twist number and page number among upward planar graphs to 5.

\begin{theorem}\label{main:5-twist}
 There is an upward planar graph $ G $ with $ \pn(G) \geq \tn(G) \geq 5 $.
\end{theorem}

\section{Bounded Width}
\label{sec:width}

Recall that the width~$w(X)$ of a subset $X \subseteq V(G)$ of the vertex set of an $st$-graph~$G$ is the largest number of vertices in~$X$ that are pairwise incomparable in~$G$.
In this section we prove that the page number is bounded by a linear function of the width.
In fact, we show a more general statement:
Given a subset~$X \subseteq V(G)$, we embed all edges of~$G[X]$ in~$\mathcal{O}(w(X))$ pages such that the book embedding of $ G[X] $ can be extended to a book embedding of $ G $, where $G[X]$ denotes the subgraph of $ G $ induced by $ X $.
In particular, the vertex ordering of $ G[X] $ is a subsequence of a topological vertex ordering of $ G $.
This generalization will be used in \cref{sec:vertices}, where we combine it with the results from \cref{sec:height}.
\cref{main:width} will follow by setting~$X = V(G)$.

The main lemma of this section (\cref{lem:width}) takes as input an~$st$-graph~$G$ and a subset~$X \subseteq V(G)$ of the vertices.
It describes how to assign all edges in~$G[X]$ to few pages.
Additionally, the lemma constructs a new~$st$-graph~$G'$ that is used in \cref{sec:vertices} to handle the remaining edges, namely those with at most one endpoint in~$X$.
We construct $ G' $ so that it is a supergraph of a subdivision of $ G $.
In particular, we have $ V(G) \subseteq V(G') $ and for every two vertices~$u,v \in V(G)$, whenever~$u \prec_G v$, then also~$u \prec_{G'} v$.
Therefore every topological ordering of~$G'$ (restricted to the vertex set of~$G$) yields a topological ordering of~$G$.
Note that for some~$u$ and~$v$, we might have~$u \prec_{G'} v$ but~$u \not\prec_{G} v$.
These additional comparabilities restrict the possible vertex orderings, making sure that the already assigned edges remain crossing-free on their respective pages, no matter which topological ordering of~$G'$ is chosen in later steps.
We remark that the subdivisions are needed to introduce comparabilities between vertices that cannot be connected by an edge in an upward planar way.
All edges in~$E(G) - E(G')$ that are removed due to the subdivisions while constructing~$G'$ are accounted for by \cref{lem:width} as well.

\begin{figure}[t]
  \centering
  \includegraphics{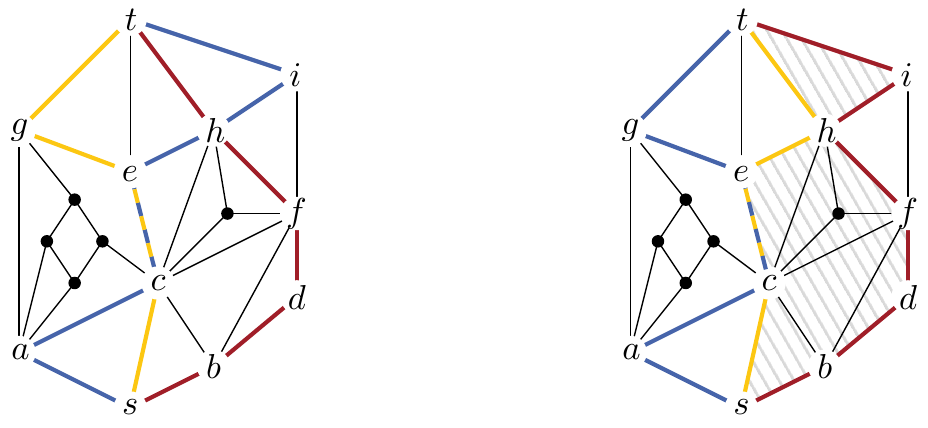}
  \caption{
    Left:~An~$st$-graph~$G$ with three paths covering the subset~$X \subseteq V(G)$ of all labeled vertices.
    Right:~The same graph with three non-crossing paths covering~$X$.
    All colored edges as well as $(a,g)$, $(e,t)$, $(c,h)$, $(b,f)$ and $(f,i)$ are intra-path-edges.
    On the other hand, $(b,c)$ and $(c,f)$ are inter-path-edges.
    The two shaded regions are the two lenses between the yellow and the red path.
  }
  \label{fig:non-crossing-path-decomposition}
\end{figure}

So consider the $st$-graph $G$ and a set $X \subseteq V(G)$.
All vertices in~$X$ can be covered by a set~$\calP$ of $st$-paths, where~$\lvert \calP \rvert = w(X)$.
To see this, consider the directed acyclic graph~$H$ with vertex set~$X$ and an edge from~$u \in X$ to~$v \in X$ if and only if~$u \prec_G v$.
Its reachability poset~$P_H$ has width~$w(P_H) = w(X)$.
By Dilworth's Theorem,~$P_H$ can be decomposed into~$w(X)$ chains, i.e.\ subsets of pairwise comparable elements.
Each of these chains can be extended to an~$st$-path in~$G$.
Given an upward planar embedding of~$G$, we can define what it means for two of these paths to cross:
Let~$P,Q \in \calP$ be two paths and~$v$ be the last vertex on the longest shared subpath beginning at~$s$ (the case~$v = s$ is possible).
Without loss of generality the next edge of~$P$ precedes the next edge of~$Q$ in the clockwise order of~$v$'s outgoing edges.
We say that~$P$ and~$Q$ \emph{cross} at another common vertex~$w$ if the next edge of~$P$ succeeds the next edge of~$Q$ in the clockwise order of~$w$'s outgoing edges.
Note that this definition allows~$P$ and~$Q$ to have common vertices and edges, even if they do not cross.
In the following we always assume~$\calP$ to be \emph{non-crossing}, meaning there is a left-to-right ordering $P_1,\ldots,P_{w(X)}$ of the $st$-paths in $\calP$ such that no two consecutive paths cross.
This assumption is justified, as a crossing between two paths~$P$ and~$Q$ at vertex~$w$ can be removed by swapping their subpaths starting at~$w$.
Thus, any set of crossing paths can be made non-crossing in every upward planar embedding (see for example the blue, yellow and red path in \cref{fig:non-crossing-path-decomposition}, which cross in the left but not in the right subfigure).

For two consecutive paths $P_i, P_{i+1} \in \calP$, a \emph{lens~$L$ between~$P_i$ and~$P_{i+1}$} is a subgraph of~$G$ enclosed by two subpaths~$P_i' \subseteq P_i$ and~$P_{i+1}' \subseteq P_{i+1}$ such that their endpoints coincide and they do not share any inner vertices.
A lens $L$ has a unique source~$s_L$ and a unique sink~$t_L$ with~$s_L, t_L \in V(P_i') \cap V(P_{i+1}')$ and~$s_L \prec_G t_L$.
As any two paths in~$\calP$ share the global source~$s$ and sink~$t$, there is at least one lens between any two consecutive paths in $\calP$.
Given the $st$-paths in $\calP$, we distinguish two kinds of edges of~$G[\calP]$, where~$G[\calP]$ is the subgraph of $G$ induced by all vertices covered by~$\calP$:
We call an edge~$e$ having both endpoints contained in the same path in~$\calP$ an \emph{intra-path-edge} ($e$ can be an edge of the path or a transitive edge).
In contrast, an \emph{inter-path-edge} $e$ has its two endpoints in two different paths in $\calP$.
We note that if~$P_i$ and~$P_{i+1}$ share some vertices and edges, it is technically possible for an edge to be an intra-path-edge and an inter-path-edge at the same time.
In this case, we consider it to be an intra-path-edge.
See \cref{fig:non-crossing-path-decomposition} for a visualization of the terminology.

\begin{lemma}
  \label{lem:width}
  Let~$G$ be an~$st$-graph and let~$X \subseteq V(G)$ be a subset of its vertices of width~$w$.
  Then there is an~$st$-graph~$G'$ with~$V(G) \subseteq V(G')$ such that:
  \begin{itemize}
    \item For every two vertices~$u,v \in V(G)$ with~$u \prec_G v$ we have~$u \prec_{G'} v$.
    \item Every topological ordering of~$G'$ admits an assignment of~$E(G'[X])$ and~$E_\Delta$ to~$14w$ pages, where $ E_\Delta = E(G) - E(G') $.
  \end{itemize}
\end{lemma}

\begin{proof}
  We start by initially setting~$G' = G$.
  As we go on, we add additional edges to~$G'$ and subdivide existing ones.
  Thus at the end~$G'$ is a supergraph of a subdivision of~$G$ and all reachabilities of~$G$ are maintained.
  All edges of~$G$ not in~$G'$ (exactly the ones that are subdivided) form the set~$E_\Delta$.

  As~$X$ has width~$w$, there is a set~$\calP = \{P_1, \ldots, P_w\}$ of non-crossing $st$-paths in~$G'$ covering all vertices of~$X$.
  Note that whenever we subdivide an edge on some path~$P \in \calP$, the new subdivision vertex and its two incident edges are added to~$P$ (replacing the subdivided edge).
  This way, the subgraph $G'[\calP]$ induced by the paths in~$\calP$ is well-defined at every step.
  Let the paths be numbered such that~$P_i$ is to the left of~$P_j$ whenever~$i < j$.
  Let~$P_i, P_{i+1} \in \calP$ be two consecutive paths in the left-to-right ordering and let~$L_1,L_2$ be two lenses between them.
  For $j = 1, 2$, let~$s_j$ and~$t_j$ denote the source, respectively the sink, of~$L_j$.
  By definition,~$s_j$ and~$t_j$ are the only vertices bounding~$L_j$ common to both~$P_i$ and~$P_{i+1}$.
  Thus we can assume without loss of generality that~$s_1 \prec_{G'} t_1 \preceq_{G'} s_2 \prec_{G'} t_2$.
  We conclude that in every topological ordering of $G'$ two edges from different lenses of $P_i,P_{i+1}$ do not cross, allowing us to deal with each lens separately and to reuse the same set of pages for all lenses between~$P_i$ and~$P_{i+1}$.

  For a single lens~$L$ between~$P_i$ and~$P_{i+1}$ we partition the inter-path-edges in~$L$ into~$\overrightarrow{E_L}$ (oriented from~$P_i$ to~$P_{i+1}$) and~$\overleftarrow{E_L}$ (oriented from~$P_{i+1}$ to~$P_i$).
  From the planarity of~$G'$ we obtain a bottom-to-top ordering $e_1,\ldots,e_\ell$ of the inter-path-edges, i.e., we order them by their endpoints along~$P_i$, using the endpoints on~$P_{i+1}$ as a tie-breaker.

  Before we actually assign the edges to pages, let us give a short overview of the strategy:
  We will consider the inter-path-edges in~$\overrightarrow{E_L}$ and~$\overleftarrow{E_L}$ separately, distributing their edges (and all edges that are subdivided in the process) to six pages each.
  This results in a total of twelve pages for all lenses between two consecutive paths~$P_i$ and~$P_{i+1}$.
  We will finish the proof by observing that each path itself (possibly with its subdivided edges) requires just two more pages.
  As there are~$w$ paths, this adds up to~$2w + 12(w-1) \leq 14w$ pages.

  In the following we only consider the inter-path-edges in~$\overrightarrow{E_L}$, the case for~$\overleftarrow{E_L}$ works symmetrically.
  Some of these edges may be transitive in~$G'[\calP]$.
  We observe that for every transitive edge~$e$ there is a non-transitive edge~$f = (v_j,w_j)$ such that~$e$ is either incident to~$v_j$ and above~$f$, or incident to~$w_j$ and below~$f$ (above and below refer to the bottom-to-top ordering of the inter-path-edges).
  See on the left of \cref{fig:edge-subdivisions} for some examples of transitive and non-transitive inter-path-edges.

  Let~$(v_1,w_1), \ldots, (v_k,w_k)$ be the subset of all inter-path-edges in~$\overrightarrow{E_L}$ that are non-transitive in~$G'[\calP]$ ordered from bottom to top.
  We observe that these edges form a matching, as otherwise at least one of them would be transitive.
  Now subdivide each~$e_j = (v_j,w_j)$ with~$j \in \{1,\ldots,k\}$ in~$G'$ and call the subdivision vertex~$u_j$.
  Further subdivide the edge of~$P_i$ outgoing from~$v_j$ and the edge of~$P_{i+1}$ incoming to~$w_j$ in~$G'$ calling the subdivision vertices~$v_j'$ and~$w_j'$, respectively.
  Note that while dealing with $\overrightarrow{E_L}$, no edge of $\overleftarrow{E_L}$ is subdivided (and vice versa).
  By upward planarity,~$v_j' \not\prec_{G'} w_j'$ and thus adding a directed path from~$w_j'$ to~$v_j'$ in~$G'$ (which we shall do next) maintains the acyclicity of $G'$.
  Additionally we ensure that~$G'$ remains a planar $st$-graph (and thus upward planar) with this new directed path (see the right of \cref{fig:edge-subdivisions}):
  Call~$E_{w,j}$ the set of edges incoming to~$w_j$ in clockwise order between~$(w_j',w_j)$ and~$(u_j,w_j)$.
  Subdivide each edge in~$E_{w,j}$ once and add a path from~$w_j'$ to~$u_j$ through the subdivision vertices in clockwise order.
  Analogously,~$E_{v,j}$ contains the edges outgoing from~$v_j$ between~$(v_j,u_j)$ and~$(v_j,v_j')$ in counterclockwise order.
  We subdivide all edges in~$E_{v,j}$ and extend the new path from~$u_j$ to~$v_j'$ through all subdivision vertices in counterclockwise order.
  Now~$w_j' \prec_{G'} v_j'$, as desired.
  Note that~$E_\Delta$ consists of all edges that were subdivided in~$G$.
  This includes all inter-path-edges and additionally some intra-path-edges and those edges incident to~$v_j$ or~$w_j$ with only one endpoint in~$\calP$.

  \begin{figure}[tb]
    \centering
    \includegraphics{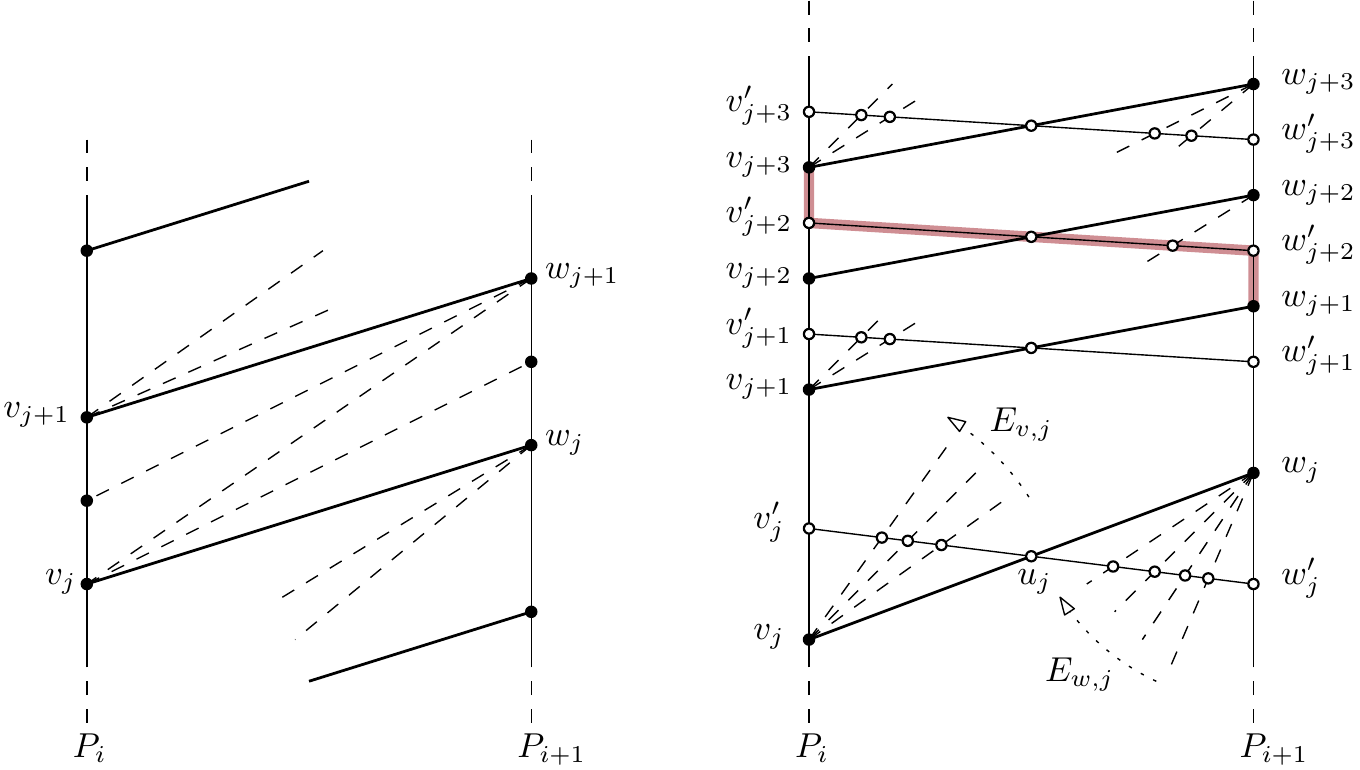}
    \caption{
      Left:~Several inter-path-edges between~$P_i$ and~$P_{i+1}$.
      Only the solid ones are non-transitive in~$G'[\calP]$.
      Right:~The comparability between~$w_j'$ and~$v_j'$ was achieved by adding a~$w_j'$-$v_j'$-path.
      To preserve that~$G'$ is a planar $st$-graph, all intersected edges are subdivided at the intersections.
      Further the comparability~$w_{j+1} \prec_{G'} v_{j+3}$ is highlighted.
      Note that the dashed edges may have only one endpoint in~$P_i$ or~$P_{i+1}$.
    }
    \label{fig:edge-subdivisions}
  \end{figure}

  We now assign the edges in~$G'[\calP]$ and~$E_\Delta$ to pages.
  Note that all edges of~$G'[\calP] \cup E_\Delta$ are inter-path-edges or intra-path-edges as both their endpoints are in $ \calP $.
  The edges in~$E_{v,j} \cup \{e_j\}$ form a star centered at~$v_j$, so they can all be assigned to the same page in any topological ordering.
  Further, all of these edges have~$v_j$ as their lower endpoint.
  In an upward planar drawing of~$G'$, all edges in~$E_{v,j}$ are inside the subregion of~$L$ enclosed by~$P_i$ and~$P_{i+1}$ to the sides and the subdivided~$e_j$ and~$e_{j+1}$ to the bottom and top.
  Thus in every topological ordering of~$G'$ they end at~$w_{j+1}$ or earlier and thus before~$v_{j+3}$ (because~$w_{j+1} \prec_{G'} w_{j+2}' \prec_{G'} v_{j+2}' \prec_{G'} v_{j+3}$, see the red path in the right of \cref{fig:edge-subdivisions}).
  Therefore the star centered at~$v_j$ can be embedded on the same page as the star centered at~$v_{j+3}$.
  Generalizing this observation, we assign~$E_{v,j} \cup \{ e_j \} $ to a page~$Q_{i,i+1}^r$ where~$r$ is the remainder of~$j$ divided by~$3$.
  With a symmetric argument all edges in~$E_{w,j}$ can be assigned to three more pages~$Q_{i+1,i}^r$.

  The intra-path-edges are left to be embedded.
  Each path~$P \in \calP$ (including the added subdivision vertices) induces a planar directed Hamiltonian graph~$H_P \subseteq G'$.
  The edges lost while subdividing can be added to~$H_P$ such that it remains planar and Hamiltonian.
  Therefore all intra-path-edges (of~$G$ and~$G'$) can be assigned to two further pages~$Q_i^1$ and~$Q_i^2$ in any topological ordering of~$G'$~\cite{BDDDMP19}.

  Let us recap, that we use twelve pages for the inter-path-edges between any two consecutive paths in~$\calP$ and two pages per path for the intra-path-edges.
  In total we get that~$2w + 12(w - 1) \leq 14w$ pages suffice for every topological ordering of~$G'$.
\end{proof}

As mentioned above, \cref{main:width} now follows as a direct corollary from \cref{lem:width} by choosing $X = V(G)$.
Let us remark that a more careful argumentation leads to a slightly better result.
We are able to show that for every~$st$-graph~$G$ we have~$\pn(G) \leq 4w(G) - 2$ by using a different strategy to embed the inter-path-edges.
However we were not able to show the more general statement of \cref{lem:width} (which we need in \cref{sec:vertices}) with this approach and hence omit this improvement of \cref{main:width} here.

\section{Bounded Height}
\label{sec:height}

In this section, we prove \cref{main:height}, which bounds the page number of any $ st $-graph in terms of its height.
Recall that the height $ h(X) $ of a subset $ X \subseteq V(G) $ of the vertex set of an $ st $-graph $ G $ is the largest number of vertices in $ X $ that are pairwise comparable in $ G $.
In combination with \cref{lem:width} from the previous section, the following lemma is central in the proof of our sublinear bound on the page number of upward planar graphs in terms of the number of vertices.
As in \cref{lem:width}, we prove a stronger statement than \cref{main:height} by considering arbitrary subsets $X$ of vertices of the graph.

\begin{lemma}\label{lem:height}
  Let~$G$ be an $st$-graph and let~$X \subseteq V(G)$ be a subset of its vertices of height~$h$.
  Then~$G$ admits a topological vertex ordering such that the size of every twist consisting of edges with at least one endpoint in~$X$ is at most~$4h$.
\end{lemma}

\begin{proof}
    Di Battista, Tamassia, and Tollis~\cite{DTT92} showed that for every $st$-graph, there is a \emph{dominance drawing}: This is a planar drawing such that between any two vertices~$u$ and~$v$, there is a path from~$u$ to~$v$ if and only if $x(u) < x(v)$ and $y(u) < y(v)$, where~$x(w)$ and~$y(w)$ denote the $x$-coordinate and $y$-coordinate of a vertex~$w$, respectively (see \cref{fig:lattice-realizer}).
    Let~$<_x$ and $<_y$ denote the vertex ordering that is given by increasing $x$-coordinates, respectively $ y $-coordinates.
    We also write $ v >_x u $ and $ v >_y u $ instead of $ u <_x v $ and $ u <_y v $, respectively.
    Most importantly, we observe that
    \begin{equation}
        u \prec_G v \quad \iff \quad u <_x v \text{ and } u <_y v.
        \label{eq:2-realizer}
    \end{equation}

    \begin{figure}
        \centering
        \includegraphics{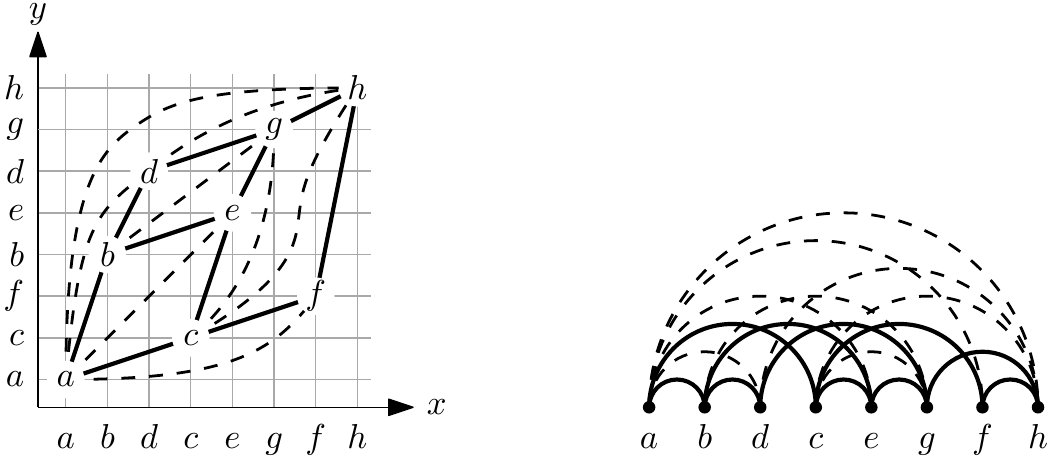}
        \caption{A dominance drawing (left) and the same graph with spine ordering $ <_x $  (right). Transitive edges are drawn dashed for better readability.}
        \label{fig:lattice-realizer}
    \end{figure}

    Now we take~$<_x$ as the linear vertex ordering for~$G$ and consider a largest twist $a_1 <_x \cdots <_x a_k <_x b_1 <_x \cdots <_x b_k$ consisting of edges in~$G$ with at least one endpoint in~$X$.
    That is, $(a_i, b_i) \in E(G)$ and we have $a_i \in X$ or $b_i \in X$ for $i=1,\ldots,k$.
    We assume for the sake of contradiction that $k > 4h$.
    By pigeonhole principle, more than~$k/2$ of the~$a_i$'s are in~$X$ or more than~$k/2$ of the~$b_i$'s are in~$X$.
    Assume the first, the latter case works symmetrically.
    In the following we refer to \cref{fig:height-contradiction}~(left), whereas the symmetric case is shown in \cref{fig:height-contradiction}~(right).
    Without loss of generality, we have $ a_1, \dots, a_{k'} \in X $, where~$k' > 2h$.
    Consider $a_1,\ldots,a_{k'}$ and their ordering with respect to~$<_y$.
    As~$k' > 2h$, by the Erd\H{o}s-Szekeres theorem there exists at least one of the following:
    \begin{itemize}
        \item A sequence $i_1 < \cdots < i_{h+1}$ of indices with $a_{i_1} <_y \cdots <_y a_{i_{h+1}}$.
        \item A sequence $i_1 < i_2 < i_3$ of indices with $a_{i_1} >_y a_{i_2} >_y a_{i_3}$.
    \end{itemize}
    The first case would give together with \eqref{eq:2-realizer} that $a_{i_1} \prec \cdots \prec a_{i_{h+1}}$, i.e.,~$h+1 $ pairwise comparable vertices in~$X$, a contradiction.
    Thus, we have the second case: Three vertices $a_{i_1},a_{i_2},a_{i_3}$ with opposing ordering with respect to~$<_x$ and~$<_y$, as illustrated in \cref{fig:height-contradiction}~(left).
    Together we have that $a_{i_2} <_x a_{i_3} <_x b_{i_1} <_x b_{i_2}$ and $a_{i_3} <_y a_{i_2} <_y a_{i_1} <_y b_{i_1}$.
    On one hand, this implies with~\eqref{eq:2-realizer} that $a_{i_3} \prec b_{i_1}$ and hence there is a path~$P$ in~$G$ from $a_{i_3}$ to $b_{i_1}$ that is monotone in $x$- and $y$\nobreakdash-coordinates, i.e.,~in the dominance drawing $P$ lies entirely inside the axis-aligned rectangle~$R$ spanned by the elements~$a_{i_3}$ and~$b_{i_1}$, see~\cref{fig:height-contradiction}~(left).
    \begin{figure}
        \centering
        \includegraphics[page = 2]{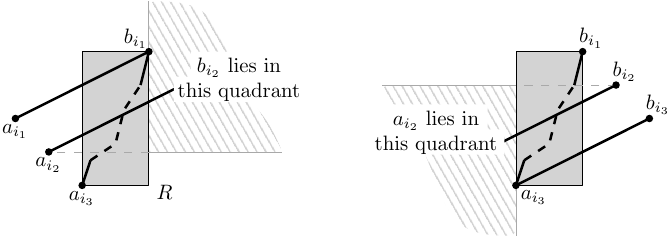}
        \caption{%
            The situation for the final contradiction in the proof of \cref{lem:height}, where
            $a_{i_1} >_y a_{i_2} >_y a_{i_3}$ (left), respectively
            $b_{i_1} >_y b_{i_2} >_y b_{i_3}$ (right, symmetric case with $ b_1, \dots, b_{k'} \in X $), by Erd\H{o}s-Szekeres.
        }
        \label{fig:height-contradiction}
    \end{figure}
    On the other hand, the edge $e = (a_{i_2}, b_{i_2})$ crosses through the rectangle~$R$ from left to right.
    Note that edge~$e$ indeed lies below~$b_{i_1}$ as it does not cross the edge $(a_{i_1}, b_{i_1})$.
    We conclude that edge~$e$ crosses path~$P$, which contradicts the planarity of the drawing.
\end{proof}

Choosing~$X = V(G)$, \cref{lem:height} gives a topological ordering of any $st$-graph~$G$ whose maximum twist size is linear in its height.
Together with \cref{obs:twist-to-page}, this proves \cref{main:height}.
We remark that for~$X = V(G)$ we can strenghten the analysis above to~$\tn(G) \leq 2h$:
As all edges have both endpoints in~$X$, we do not need to apply the pigeonhole principle to get that at least half the twisting edges have their lower (or equally their upper) endpoint in~$X$, thus saving a factor of~$2$.

\section{Bound in Terms of the Number of Vertices}
\label{sec:vertices}

In this section we combine our approaches of bounding the page number in terms of width and height and obtain the first sublinear upper bound on the page number of upward planar graphs and planar posets.
We prove \cref{main:vertices} which states that the page number of~$n$\nobreakdash-vertex upward planar graphs is $\mathcal{O}(n^{2/3} \log^{2/3}(n))$.

\begin{proof}[Proof of \cref{main:vertices}]
    Let~$G$ be an~$n$-vertex upward planar graph.
    Without loss of generality we may assume that~$G$ is an~$st$-graph~\cite{Kel87}.
    We first identify vertices that can be covered by few long directed paths and use \cref{lem:width} to embed the subgraph induced by these paths.
    We then apply \cref{lem:height} to the remaining vertices to find a topological ordering that admits an assignment of the remaining edges to few pages.
    However, as \cref{lem:width} introduces new directed reachabilities to the graph, we have to pick the first vertex set in a sequential way.

    We construct a sequence $ G_0, G_1, \dots $ of graphs and a sequence $ L_0 \subseteq  L_1 \subseteq \dots $ of sets containing the vertices of \enquote{long} directed paths in the respective graphs, starting with~$G_0 = G$ and~$L_0 = \emptyset$.
    We thereby ensure that $ V(G_i) \subseteq V(G_{i + 1}) $ and that the width of~$L_i$ in~$G_i$ is at most~$i$ for each~$i \geq 0$.
    Let $ \ell = n^{2/3} / \log^{1/3}(n) $; we use this threshold to decide which paths are considered \emph{long paths}.
    For ease of notation, let $ E_\Delta(i, i + 1) = E(G_i) - E(G_{i + 1}) $ denote the set of edges of~$G_i$ that is removed when defining the next graph~$G_{i + 1}$.
    We write~$G[X]$ for the subgraph of $ G $ that is induced by~$X \cap V(G)$, where~$X$ is a set of vertices of~$G_i$ which may include vertices that are not in~$G$.

    Assume that~$G_i$ and~$L_i$ are already defined and that there is an $st$\nobreakdash-path~$P$ in~$G_i$ that contains at least~$\ell$ vertices of~$G$ that are not contained in~$L_i$.
    We include the vertices of~$P$ in the next set~$L_{i + 1}$.
    That is, we define $ L_{i + 1} = L_i \cup V(P) $.
    Note that adding the vertex set of a directed path to a set of vertices increases the width by at most~$1$.
    Hence, the width of~$L_{i + 1}$ in~$G_i$ is at most~$i + 1$.
    Now apply \cref{lem:width} to~$G_i$ and~$L_{i + 1}$ and obtain an $st$-graph~$G_{i + 1}$ with $V(G_i) \subseteq V(G_{i+1})$.
    By \cref{lem:width}, every topological ordering of~$G_{i + 1}$ admits an assignment of $ E(G_{i + 1}[L_{i + 1}]) \cup E_\Delta(i, i + 1) $ to~$14(i + 1)$ pages.
    As~$L_{i + 1}$ can be covered by~$i + 1$ paths in~$G_i$, the same holds in~$G_{i + 1}$ as the reachabilities are preserved.
    Thus, the width of~$L_{i + 1}$ in~$G_{i + 1}$ is at most~$i + 1$.
    We remark that \cref{lem:width} does not give any guarantees on the size of $ G_{i + 1} $.
    However, the number of pages the lemma uses is bounded in terms of the width, which increases only by 1 in each step.

    Let~$t$ denote the largest~$i$ for which~$G_i$ and~$L_i$ are defined, i.e., there is no path in~$G_t$ that contains at least~$\ell$ vertices of the initial graph~$G$ that are not covered by~$L_t$.
    Note that $t \leq n / \ell = n^{1/3}\log^{1/3}(n)$, because we add at least~$\ell$ vertices of~$G$ in each round.
    We claim that every topological ordering of~$G_t$ restricted to~$V(G)$ admits an assignment of the edges of~$G[L_t]$ to~$\mathcal{O}(t^2)$ pages.
    To this end, fix an arbitrary topological ordering~$<_t$ of~$G_t$ and consider the restriction~$<$ of~$<_t$ to the vertex set of~$G$.
    Observe that~$<$ is a topological ordering of~$G$ as directed reachabilities in~$G$ are maintained in~$G_t$.
    For $ i = 0, \dots, t - 1 $, let $ \calQ_{i, i + 1} $ denote the set of~$14(i + 1)$ pages used by \cref{lem:width} when applied to~$G_i$.
    We restrict the pages to contain only edges of~$G$.
    Observe that the edges in $ E_\Delta(i, i + 1) \cap E(G) $ are embedded in some page of~$\calQ_{i, i + 1}$ for each $ i = 0, \dots, t - 1 $.
    Now, let $ E_t $ denote the remaining edges of $ G[L_t] $.
    Note that these edges are contained in $ G_t[L_t] $ and thus are embedded in some page of $ \calQ_{t - 1, t} $ by \cref{lem:width}.
    We conclude that the union of all $ \calQ_{i, i + 1} $ covers all edges of $ G[L_t] $ with
    $ \sum_{i = 0}^{t - 1} |\calQ_{i, i + 1}| = \sum_{i = 0}^{t - 1} 14 (i + 1) = 7 t (t + 1) $ pages.

    It is left to embed the set~$E_S$ of edges in~$G$ that are also contained in~$G_t$ and have at most one endpoint in~$L_t$, i.e., at least one endpoint in $S = V(G)-L_t$.
    Note that there are two types of edges in~$G$ having exactly one endpoint in $ L_t $: First, edges that are subdivided by \cref{lem:width} and thus are not contained in $ G_t $. 
    These edges are already embedded by \cref{lem:width}.
    And second, edges that are not touched by \cref{lem:width} and thus belong to $ G_t $.
    These edges, together with edges having no endpoint in $ L_t $, are handled by \cref{lem:height} as follows.
    Recall that there is no path in~$G_t$ with at least~$\ell$ vertices of~$G$ that are not contained in~$L_t$, i.e., the height of~$S$ in~$G_t$ is less than~$\ell$.
    Applying \cref{lem:height} to~$S$ and~$G_t$ yields a topological ordering~$<_t$ of~$G_t$ such that the edges in~$E_S$ form twists of size at most~$4 \ell$.
    By \cref{obs:twist-to-page}, the same vertex ordering admits an assignment of the edges in~$E_S$ to~$\mathcal{O}(\ell \log(\ell)) = \mathcal{O}(\ell \log(n)) $ pages.
    Restricting~$<_t$ to~$G$ and combining the page assignment of~$E_S$ with the page assignment of~$G[L_t]$, we obtain a book embedding of~$G$ with $\mathcal{O}(\ell \log(n) + t^2) = \mathcal{O}(n^{2/3} \log^{2/3}(n))$ pages (recall that~$\ell = n^{2/3} / \log^{1/3}(n)$ and $t \leq n / \ell = n^{1/3} \log^{1/3}(n)$).
\end{proof}

In view of \cref{lem:height} which bounds the twist number instead of the page number, the question arises whether our bound in terms of the number of vertices can be decreased by improving this step.
We also remark that by choosing $ \ell = n^{2/3} $, we obtain that every upward planar graph admits a topological vertex ordering whose maximum twist size is $ \mathcal{O}(n^{2/3}) $.
That is, any improvement in bounding the page number of upward planar graphs in terms of their twist number also improves our result.
Such an improvement, however, needs to make use of the structure of the graph and the constructed vertex ordering as Davies' result is asymptotically tight.

\section{Lower Bound}
Recall that the twist number $\tn(G)$ of a directed acyclic graph $G$ is the maximum $k$ for which every topological ordering of $G$ contains $k$ pairwise crossing edges.
In this section, we construct an upward planar graph whose twist number, and therefore in particular its page number, is at least~$5$.
This improves on the previously best known bound of an upward planar graph that requires four pages (but has twist number 3) by Hung~\cite{Hun93}.
We remark that the second author~\cite{Mer20} improved on our upward planar graph by transforming it into a planar poset whose twist number and page number are at least~$5$.

\begin{figure}
  \centering
  \includegraphics{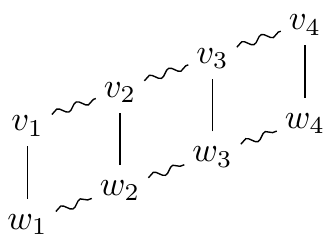}
  \hspace{4em}
  \includegraphics{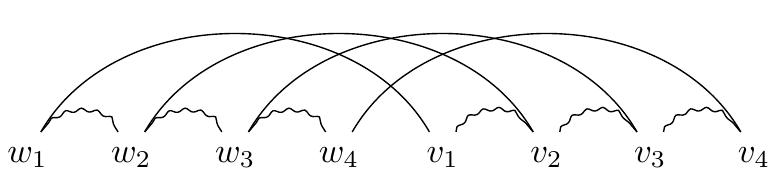}
  \caption{
    A $4$-fence from~$v_1$ to~$w_4$ and a topological ordering with~$w_4 < v_1$ yielding a $4$-twist.
    Note that $(v_1, w_1), \dots, (v_4, w_4)$ are edges, while all other shown reachabilites may be due to paths.
  }
  \label{fig:fence}
  \label{fig:fence-twist}
\end{figure}

We identify a structure that can lead to large twists if the spine ordering is not chosen carefully.
By adding additional edges, any topological ordering of the augmented graph avoids these twists.
For~$k \geq 2$, a \emph{$k$-fence (from~$v_1$ to~$w_k$)} consists of~$2k$ distinct vertices $v_1 \prec \dots \prec v_k$ and $w_1 \prec \dots \prec w_k$ together with the edges $(w_i, v_i)$ for each $i = 1, \dots, k$.
The edges $(w_i,v_i)$ are called \emph{fence edges}.
If~$k$ is not important, we simply say \emph{fence}.
\Cref{fig:fence}~(left) shows a $4$-fence.
Observe that~$v_1$ and~$w_k$ are not necessarily comparable.
However, we show that~$v_1$ must preceed~$w_k$ in every spine ordering that has no~$k$-twist.
By transitivity, every~$v \preceq v_1$ must therefore also preceed every~$w \succeq w_k$.

\begin{observation}
  \label{obs:fence}
  Every topological ordering $ < $ of a~$k$-fence from~$v_1$ to~$w_k$ in which~$w_k < v_1$ has a~$k$-twist.
\end{observation}

\begin{proof}
  Assuming $w_k < v_1$, we obtain $w_1 < \dots < w_k < v_1 < \dots < v_k$ as the unique topological ordering.
  Hence, the fence edges form a~$k$-twist.
  See \cref{fig:fence-twist} (right) for an illustration.
\end{proof}

Given an upward planar graph~$G$, we augment it with additional edges that indicate how to avoid the~$k$-twists that otherwise are present in~$k$-fences.
Let~$k \geq 2$ and consider a~$k$-fence from~$v_1$ to~$w_k$ in~$G$ such that $v_1 \not\prec w_k$.
We add a new edge~$(v_1,w_k)$ forcing~$v_1 < w_k$, as every topological ordering with $w_k < v_1$ yields a~$k$-twist (\cref{obs:fence}).
Since adding edges increases the set of reachabilities in~$G$, new fences might emerge with each newly added edge.
Here we consider only fences for which the fence edges~$(v_i,w_i)$ for~$i = 1,\ldots, k$ are edges of the original graph~$G$, whereas the reachabilities along the two paths~$v_1 \prec \dots \prec v_k$ and~$w_1 \prec \dots \prec w_k$ might consist (partly) of new edges.
We continue adding new edges to each current and future~$k$-fence.
This process terminates, as there is only a finite number of possible comparabilities between the unchanged number of vertices.
The resulting graph is denoted by~$G_k^\ast$ and contains no $k$-fence from $v_1$ to $w_k$ with $v_1 \not\prec w_k$.

\begin{figure}
  \centering
  \includegraphics{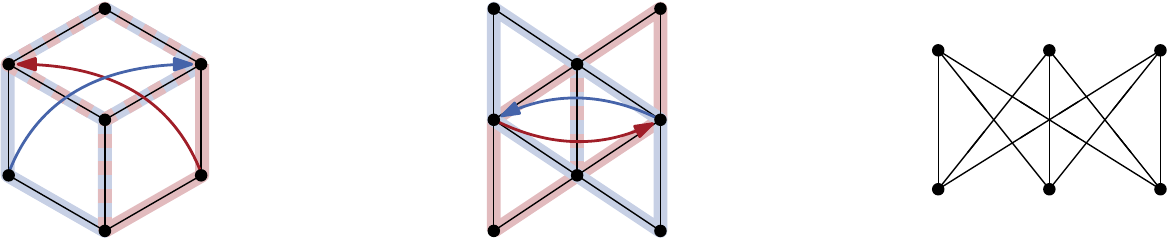}
  \caption{
    Examples of $G_3^\ast$.
    The edges in $E(G_3^\ast)-E(G)$ are drawn thick and the $3$-fences are highlighted.
    Left:~Taking a topological ordering of $G_3^\ast$ shows that $\tn(G)\leq 2$.
    Middle:~$G_3^\ast$ is not acyclic and hence $\tn(G)> 2$.
    Right:~$G_3^\ast = G$ as there is no $3$-fence, but still $\tn(G)>2$.
  }
  \label{fig:fence-examples}
\end{figure}

Let us refer to \cref{fig:fence-examples} for some illustrative examples.
Note that even if $G$ is upward planar, then~$G_k^\ast$ is not necessarily upward planar; possibly not planar, nor even acyclic.
We emphasize that the new edges in $G_k^\ast$ are not part of~$E(G)$, and as such, need not be assigned to any page in a book embedding.
Their sole purpose is to restrict the set of possible topological orderings of $G$ to those of $G_k^\ast$.

\Cref{obs:fence} shows that a topological ordering of~$G$ which is not a topological ordering of~$G_{k + 1}^\ast$ yields a $(k + 1)$-twist.
In particular, $G_{k + 1}^\ast$ being acyclic is a necessary condition for $G$ admitting a~$k$\nobreakdash-page book embedding.

\begin{corollary}
  \label{lem:G*-top}
  Every book embedding of a directed acylic graph~$G$ without a $(k + 1)$-twist (in particular every $k$-page book embedding) uses a topological ordering of $G_{k + 1}^\ast$ as spine ordering.
\end{corollary}

However, using a topological ordering of~$G_k^\ast$ as a spine ordering is not sufficient to avoid $k$-twists; see e.g., the right of \cref{fig:fence-examples}.
Quite the contrary, we find that for some small~$k$, the augmented graph~$G_k^\ast$ might be cyclic and therefore not have any topological ordering at all.
And even if~$G_k^\ast$ is acyclic, choosing any topological ordering of~$G_k^\ast$ can inescapably lead to arbitrarily large twists (which are not due to fences) even if the graph admits a book embedding with few (but more than~$k$) pages.
We shall force such a situation in our construction of an upward planar graph with twist number at least $5$, which then proves \cref{main:5-twist}.

\begin{figure}[tb]
  \centering
  \includegraphics{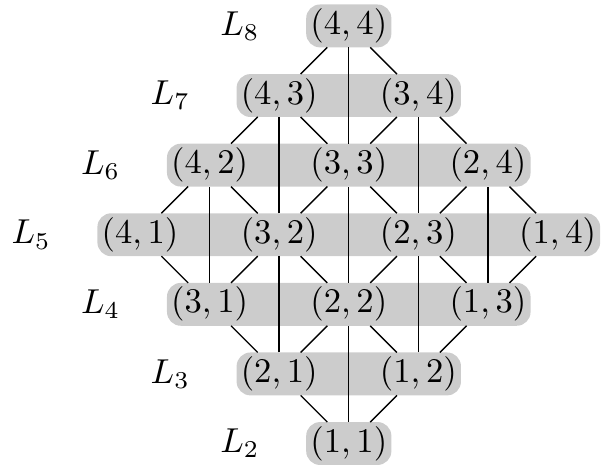}
  \caption{
    A~$4 \times 4$ upward grid with levels $ L_2, \dots, L_8 $.
    Consider the vertex~$(2, 2)$.
    The first left upper vertex is~$(3, 2)$, the second left upper vertex is~$(4, 1)$, the first right upper vertex is~$(2, 3)$, and the second right upper vertex is~$(1, 4)$.
  }
  \label{fig:upward-grid}
\end{figure}

For any integer~$n > 0$, we define an \emph{$n \times n$ upward grid}~$\Grid_n$ as follows (see \cref{fig:upward-grid}).
The vertices of~$\Grid_n$ are the tuples~$(\ell, r)$ of integers with $1 \leq \ell, r \leq n$.
The vertices are partitioned into \emph{levels}, where level~$L_h$ contains the vertices~$(\ell, r)$ with~$\ell + r = h$.
The edge set of~$\Grid_n$ consists of three subsets.
There are \emph{left edges} of the form $\gedge{(\ell,r)}{(\ell+1,r)}$ for each $r = 1, \dots, n$ and $\ell = 1, \dots, n - 1$.
Symmetrically, the edges $\gedge{(\ell,r)}{(\ell,r+1)}$ for $\ell = 1, \dots, n$ and $r = 1, \dots, n - 1$ are called \emph{right edges}.
Finally, we have edges $\gedge{(\ell,r)}{(\ell+1,r+1)}$ for $1 \leq \ell, r \leq n - 1$ and call them \emph{vertical edges}.

Consider a vertex~$v = (\ell_v, r_v)$ in some level~$L_h$ of an upward grid.
A vertex~$w = (\ell_w, r_w)$ in level~$L_{h + 1}$ is called an \emph{$i$-th left (right) upper vertex of~$v$} if $\ell_w = \ell_v + i$ ($r_w = r_v + i$).
A vertex that is an $i$-th left upper vertex or an $i$-th right upper vertex of~$v$ is also called an~\emph{$i$-th upper vertex of~$v$}.
Note that every vertex in~$L_{h + 1}$ is an~$i$-th upper vertex of~$v$ for some~$i > 0$.

\begin{figure}[tb]
  \centering
  \includegraphics[height = 26ex]{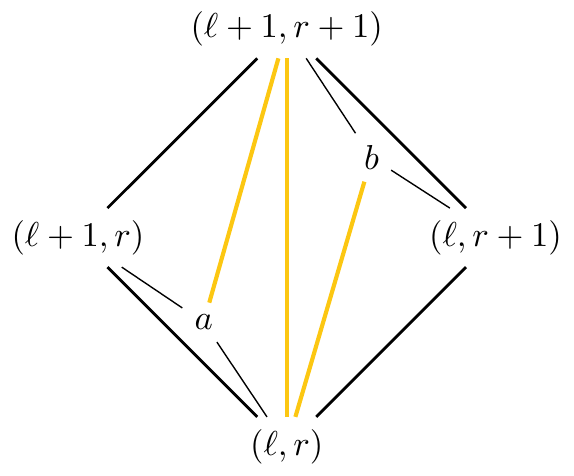}
  \hfill
  \includegraphics[height = 26ex]{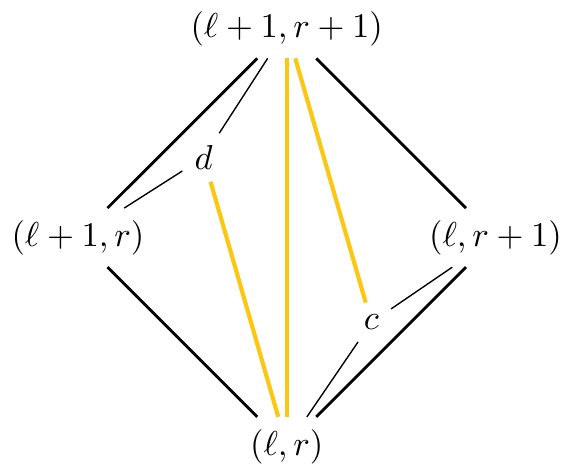}
  \hfill
  \includegraphics[height = 26ex]{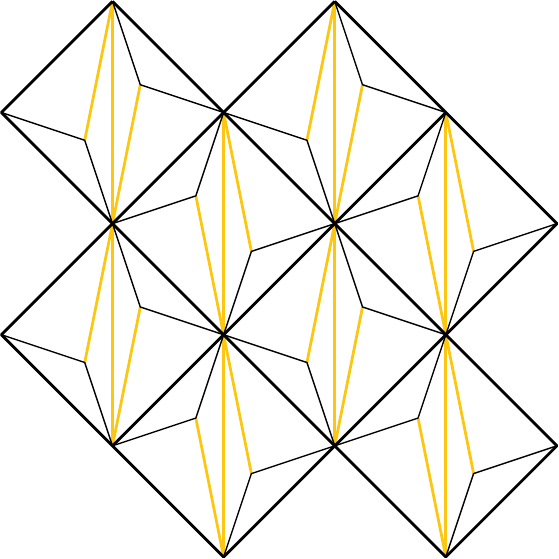}
  \caption{
    Parts of an N-grid with N-vertices $a = a_{\ell, r}$, $b = b_{\ell, r}$, $c = c_{\ell, r}$, and $d = d_{\ell, r}$, where $ \ell - r $ is even (left), respectively odd (middle).
  }
  \label{fig:N-cell}
  \label{fig:mirrored-N-cell}
  \label{fig:N-grid}
\end{figure}

Based on an $n \times n$ upward grid, we define an $n \times n$ N-grid, which we denote by~$N_n$, for any integer~$n>0$.
We shall show in this section that every $n \times n$ N-grid has a $5$-twist in every topological ordering, provided $n$ is large enough.
The~$n \times n$ N-grid~$N_n$ contains an~$n \times n$ upward grid~$\Grid_n$ as an induced subgraph and an additional vertex in each inner face of~$\Grid_n$.
The additional vertices are called \emph{N-vertices}, whereas the vertices that belong to~$\Grid_n$ are called \emph{grid vertices}.
See \cref{fig:N-grid} for an illustration.
Consider two triangles in~$\Grid_n$ that share a vertical edge.
That is, they consist of vertices~$(\ell, r)$, $(\ell + 1, r)$, $(\ell, r + 1)$, and $(\ell + 1, r + 1)$ as shown in \cref{fig:N-cell,fig:mirrored-N-cell}.
If $\ell - r$ is even, then we insert a vertex~$a=a_{\ell, r}$ into the left triangle and add edges $\gedge{(\ell, r)}{a}$, $\gedge{a}{(\ell + 1, r)}$, and $\gedge{a}{(\ell + 1, r + 1)}$.
In addition, we insert a vertex~$b=b_{\ell, r}$ together with the edges $\gedge{(\ell, r)}{b}$, $\gedge{(\ell, r + 1)}{b}$, and $\gedge{b}{(\ell + 1, r + 1)}$ into the right triangle in this case.
If~$\ell - r$ is odd, then we insert vertices~$c=c_{\ell, r}$ and~$d=d_{\ell, r}$ into the right, respectively left, triangle and add edges
$\gedge{(\ell, r)}{c}$, $\gedge{c}{(\ell, r + 1)}$, $\gedge{c}{(\ell + 1, r + 1)}$, $\gedge{(\ell, r)}{d}$, $\gedge{(\ell + 1, r)}{d}$, and $\gedge{d}{(\ell + 1, r + 1)}$.
The definitions of levels and upper vertices remain as in~$\Grid_n$, the N-vertices do not belong to any level.
Observe that every N-grid is upward planar.
Whenever we refer to an embedding of an N-grid, we assume the upward grid induced by the grid vertices to be embedded in the canonical way shown in \cref{fig:upward-grid} and the N-vertices to be placed in the respective triangular inner faces as shown in \cref{fig:N-grid}.

The rest of this section is devoted to proving that every topological ordering of a sufficiently large N-grid yields a $5$-twist.
For this, we consider the graph~$N_{n, 5}^\ast$ that results from augmenting~$N_n$ via $5$-fences as described above.
By \cref{lem:G*-top}, every topological ordering of~$N_n$ that is not a topological ordering of $N_{n, 5}^\ast$ yields a $5$-twist.
Hence, we only need to consider topological orderings of $N_{n, 5}^\ast$.
We say that two levels~$L_i,L_j$ ($2 \leq i < j \leq 2n$) are \emph{separated} by a topological ordering~$<$, if for all grid vertices~$(\ell_i,r_i) \in L_i$ and~$(\ell_j,r_j) \in L_j$ we have~$(\ell_i,r_i) < (\ell_j,r_j)$.
We write~$L_i < L_j$ in this case.
We call a topological ordering~$ < $ of an N-grid \emph{level-separating} if it separates every two consecutive levels, i.e., we have $ L_2 < \dots < L_{2n} $.
We also say that $ < $ \emph{separates the levels of the N-grid} in this case.
The next \lcnamecref{lem:separating-levels-N} shows that we can assume the levels of~$N_n$ to be separated if the vertex ordering is a topological ordering of $N_{n, 5}^\ast$.

\begin{lemma}
  \label{lem:separating-levels-N}
  For every~$n > 0$, there is an~$n' \geq n$ such that for every topological ordering~$<$ of~$N_{n', 5}^\ast$ we find a copy of $N_n \subseteq N_{n'} $ whose levels are separated by~$<$.
\end{lemma}

\begin{proof}
  We choose~$n' = n + 2 (n - 1)$ and use induction on~$i = 1, \ldots, n$.
  For each $ i $, we identify a set of vertices $ V_i \subseteq V(N_{n', 5}^\ast) $ such that each grid vertex in $ V_i $ has an outgoing edge to all its~$j$-th upper vertices that are contained in $ V_i $ for each~$j = 1, \ldots, i$.
  We thereby ensure $ V_i \subseteq V_{i - 1} $ for $ i > 1 $ and that $ V_i $ induces a copy of $ N_{n + 2(n - i)} $ in $ N_{n'} $.
  Finally we show that in $ N_{n', 5}^\ast[V_n] $ every grid vertex reaches every vertex of $ V_n $ in the subsequent level, and thus $ V_n $ induces the desired copy of $ N_n $ in~$ N_{n'} $.

  For~$i = 1$, define~$ V_1 = V(N_{n',5}^\ast) $.
  Observe that in each N-grid, every grid vertex is adjacent to its first left upper vertex via a left edge and to its first right upper vertex via a right edge, which settles the base case.
  Now let~$i > 1$ and assume that all grid vertices in $ V_{i - 1} $ reach all~$j$-th upper vertices also contained in $ V_{i - 1} $ for each~$j \leq i - 1$.
  Consider the subgraph $ N_{n'} $ of $ N_{n', 5}^\ast $ on the same vertex set but without the augmenting edges.
  To obtain~$ V_i $ from $ V_{i - 1} $, we drop all grid vertices incident to the outer face of~$ N_{n'}[V_{i - 1}] $ and then remove all N-vertices that are now incident to the outer face.
  Refer to \cref{fig:lower-bound-induction} to see how~$V_i$ lies in~$V_{i-1}$.
  \begin{figure}
    \centering
    \includegraphics{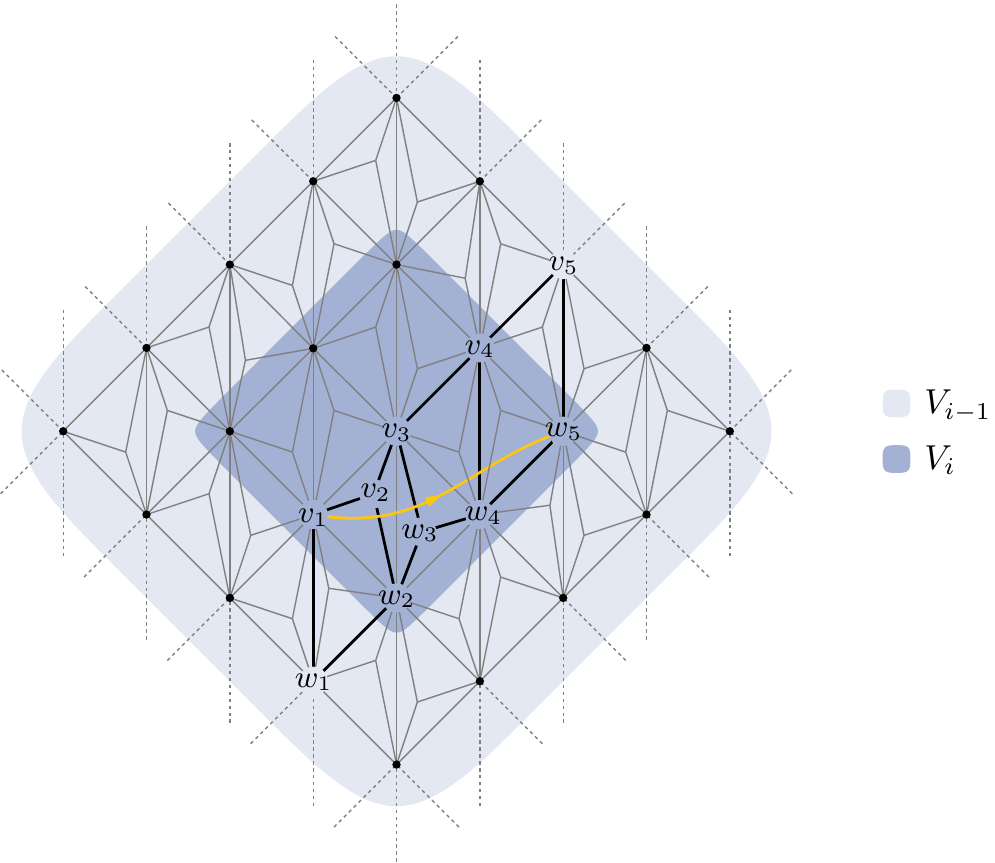}
    \caption{
      An inner N-grid~$N_{n'}[V_i]$ (darkblue) inside an outer N-grid~$N_{n'}[V_{i-1}]$ (lightblue).
      Observe that the shown~$5$-fence has vertices~$w_1$ to~$v_5$ outside~$N_{n'}[V_i]$, but the yellow edge is inside.
    }
    \label{fig:lower-bound-induction}
  \end{figure}
  Note that every grid vertex in~$N_{n'}[V_i]$ has an incoming vertical edge and an outgoing vertical edge in~$ N_{n'}[V_{i - 1}] $.
  Also observe that $ V_i $ induces an N-grid whose size is reduced by 2 in both directions compared to the N-grid $ N_{n'}[V_{i - 1}] $.

  We next find a~$5$-fence from each grid vertex of~$ V_i $ to its~$i$-th upper vertices in~$ V_i $.
  Consider a grid vertex~$v = (\ell_v, r_v) \in V_i $ .
  Without loss of generality, we assume that~$\ell_v - r_v$ is even.
  Swap left and right otherwise.
  Let $ w = (\ell_w, r_w) \in V_i $ denote the~$i$-th right upper vertex of~$v$ (if it exists).
  By definition of an~$i$-th right upper vertex, we have~$r_w = r_v + i$.
  As the two vertices are in consecutive levels, we have~$\ell_v + r_v = h$ and~$ \ell_w + r_w = h + 1$, where~$L_h$ is the level of~$(\ell_v, r_v)$.
  It follows that~$\ell_w = \ell_v - i + 1$.

  Now, consider the vertices
  \begin{align*}
      w_1 & = (\ell_v - 1, r_v - 1), \\
      w_2 & = (\ell_v - 1, r_v), \\
      w_3 & = c_{\ell_v - 1, r_v}, \\
      w_4 & = (\ell_v - 1 , r_v + 1), \text{ and} \\
      w_5 & = (\ell_v - i + 1, r_v + i) = w.
  \end{align*}%
  \begin{figure}%
      \centering%
      \includegraphics{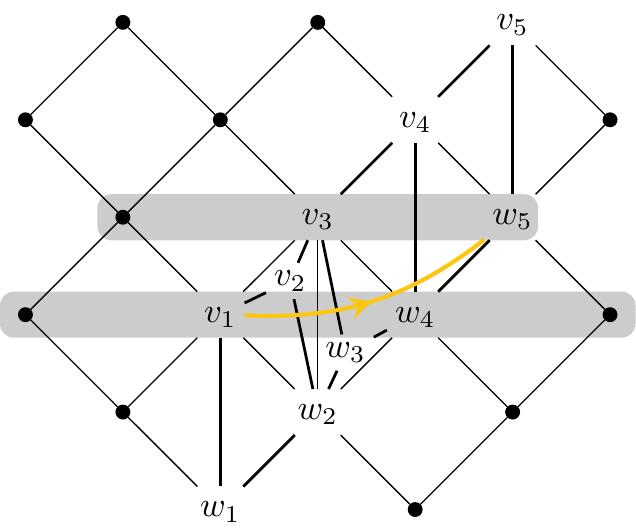}%
      \hfill%
      \includegraphics{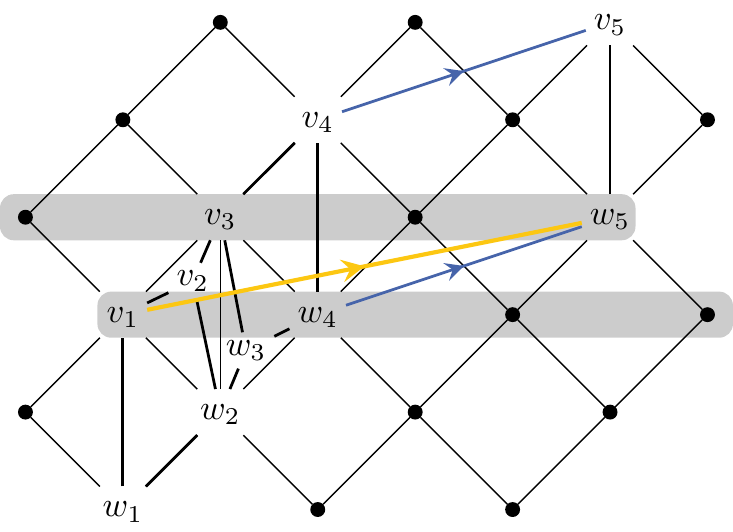}%
      \caption{
        A $5$-fence from~$v = v_1$ to~$w = w_5$, where~$w$ is the second/third right upper vertex of~$v$.
        The blue edges $ (v_4, v_5) $ and $ (w_4, w_5) $ exist by induction.
      }%
      \label{fig:right-upper-vertices}%
  \end{figure}%
  See \cref{fig:right-upper-vertices} for an illustration.
  These five vertices form the lower part of the desired $5$\nobreakdash-fence.
  Note that~$w_1$ is not necessarily in $ V_i $ but is connected to~$v$ by a vertical edge in~$N_{n'}[V_{i-1}]$ and thus is contained in $ V_{i - 1} $ (see again \cref{fig:lower-bound-induction}, where $ v = v_1 $).
  We next observe that~$w_1, \dots, w_5$ are pairwise comparable.
  The first four vertices induce a path in $ N_{n'} $.
  The edge~$(w_4,w_5)$ exists in~$ N_{n', 5}^\ast[V_{i - 1}] $ by the induction hypothesis since~$w_5$ is an~$(i-1)$-th upper vertex of~$w_4$.
  To see this, observe that~$w_4$ and~$w_5$ are in consecutive levels as~$(\ell_v - i + 1 + r_v + i) - (\ell_v - 1 + r_v + 1) = 1$ and their~$r$-coordinates differ by exactly~$i - 1$.

  Now, consider the vertices
  \begin{align*}
      v_1 & = (\ell_v, r_v) = v, \\
      v_2 & = d_{\ell_v - 1, r_v}, \\
      v_3 & = (\ell_v, r_v + 1), \\
      v_4 & = (\ell_v, r_v + 2), \text{ and} \\
      v_5 & = (\ell_w + 1, r_w + 1) = (\ell_v - i + 2, r_v + i + 1).
  \end{align*}
  These five vertices serve as the upper part of the 5-fence from $ v $ to $ w $.
  Again, we find that there is a path connecting the five vertices in $ N_{n', 5}^\ast[V_{i - 1}] $.
  First, the edges~$(v_1,v_2)$ and~$(v_2,v_3)$ exist by construction of an N-grid.
  The edge~$(v_3,v_4)$ is a right edge in~$N_{n'}$.
  We again remark that~$v_5$ is not necessarily in $ V_i $ but is connected to~$w$ via a vertical edge and thus is contained in $ V_{i - 1} $.
  We obtain the remaining edge~$(v_4,v_5)$ by induction as $ v_5 $ is an $ (i - 1) $-th right upper vertex of $ v_4 $.
  Thus, we find a 5-fence from~$v$ to~$w$ using $(v_1,w_1), \dots, (v_5,w_5)$ as fence edges.

  \begin{figure}%
      \centering%
      \includegraphics{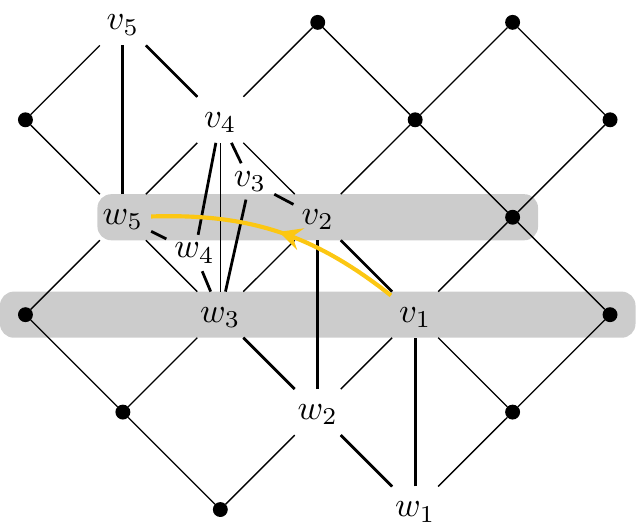}%
      \hfill%
      \includegraphics{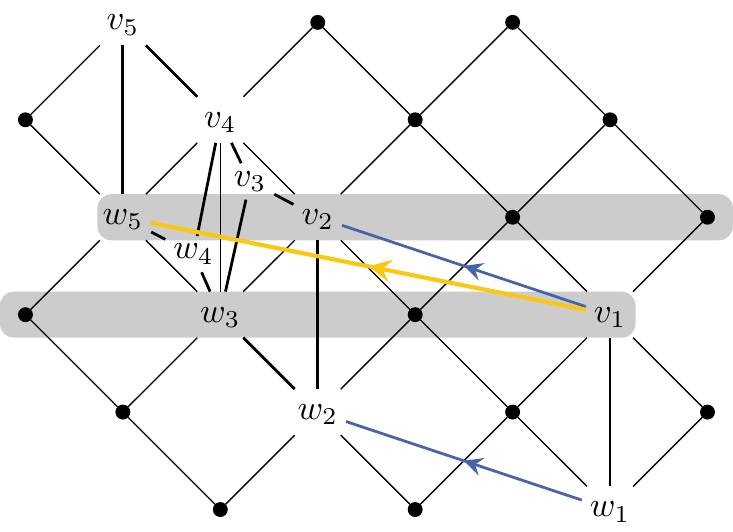}%
      \caption{%
        A 5-fence from~$v = v_1$ to~$w = w_5$, where~$w$ is the second/third left upper vertex of~$v$.
        The blue edges $ (v_1, v_2) $ and $ (w_1, w_2) $ exist by induction.
      }%
      \label{fig:left-upper-vertices}%
  \end{figure}

  The proof for the~$i$-th left upper vertex works nearly symmetrically.
  In contrast to the right upper vertex, we first use the edges obtained by the induction hypothesis and then the edges of $ N_{n'} $ to find the two paths of the 5-fence.
  Let~$w = (\ell_w, r_w) = (\ell_v + i, r_v - i + 1) \in V_i $ denote the~$i$-th left upper vertex of~$v$ (if it exists).
  We find a~$5$-fence from~$v$ to~$w$ using the vertices
  \begin{align*}
      w_1 & = (\ell_w - i - 1, r_w + i - 2) = (\ell_v - 1, r_v - 1), \\
      w_2 & = (\ell_w - 2, r_w), \\
      w_3 & = (\ell_w - 1, r_w), \\
      w_4 & = a_{\ell_w - 1, r_w}, \text{ and} \\
      w_5 & = (\ell_w, r_w) = w \\
  \intertext{for the lower part, while the upper part is formed by the vertices}
      v_1 & = (\ell_w - i, r_w + i - 1) = v, \\
      v_2 & = (\ell_w - 1, r_w + 1), \\
      v_3 & = b_{\ell_w - 1, r_w}, \\
      v_4 & = (\ell_w, r_w + 1), \text{ and} \\
      v_5 & = (\ell_w + 1, r_w + 1).
  \end{align*}
  We refer to \cref{fig:left-upper-vertices} for an illustration.
  Note that the coordinates of~$w_3$ have an even difference as~$(\ell_w - 1) - r_w = (\ell_v + i - 1) - (r_v - i + 1) = \ell_v - r_v + 2i - 2 $, which means that the claimed~$a$- and~$b$-vertices indeed exist.
  The edges~$(w_1,w_2)$ and~$(v_1,v_2)$ exist by induction as their upper endpoints are $(i-1)$\nobreakdash-th left upper vertices of the lower endpoints.
  The other vertices are connected by two paths using only edges of~$N_{n'}$.
  We again obtain a~$5$-fence using the edges $(v_1,w_1), \dots, (v_5,w_5)$ as fence edges.

  To conclude the proof, recall that $ V_n $ induces a copy of $ N_n $ in $ N_{n'} $.
  Observe that in $ N_n $, no vertex has an $ i $-th upper vertex for $ i > n $.
  Thus by the induction above, we have that in $ N_{n', 5}^\ast[V_n] $ every grid vertex reaches all its upper vertices that are contained in $ V_n $, i.e., all vertices of the subsequent level of $ N_n $.
  Therefore, the levels of $ N_{n', 5}^\ast[V_n] $ are separated by every topological ordering of $ N_{n', 5}^\ast $.
\end{proof}




Having \cref{lem:separating-levels-N}, we know that we may assume the levels of an N-grid $ N_n $ to be separated when we try to avoid 5-twists, i.e., when we consider topological orderings of $ N_{n, 5}^\ast $.
The next lemma, however, shows that separated levels imply not only 5-twists but arbitrarily large twists, finishing the proof of \cref{main:5-twist}.

\begin{lemma}
  \label{lem:G5-ordering}
  For every~$p \geq 0$, there is an~$n$ such that every level-separating topological ordering~$<$ of~$N_n$ yields a $(p + 1)$\nobreakdash-twist.
  In particular, $<$ does not admit a~$p$\nobreakdash-page book embedding.
\end{lemma}

\begin{proof}
  Let~$r = p^3 + 1$ and~$n = r + 1$.
  We identify~$r$ triangles in~$N_n$, each of which has exactly one vertex in each of the three levels~$L_n$, $L_{n+1}$ and~$L_{n+2}$.
  Observe that each of these levels has at least~$r$ vertices.
  For $i = 1, \dots, r$, we define the triangle~$T_i$ consisting of the vertices $x_i = (n - i, i) \in L_n$, $y_i = (n - i, i + 1) \in L_{n + 1}$, and $z_i = (n - i + 1, i + 1) \in L_{n + 2}$.
  See \cref{fig:triangles-in-grid} for an example.
  By our assumption we have~$L_n < L_{n+1} < L_{n+2}$.

  \begin{figure}
    \centering
    \includegraphics{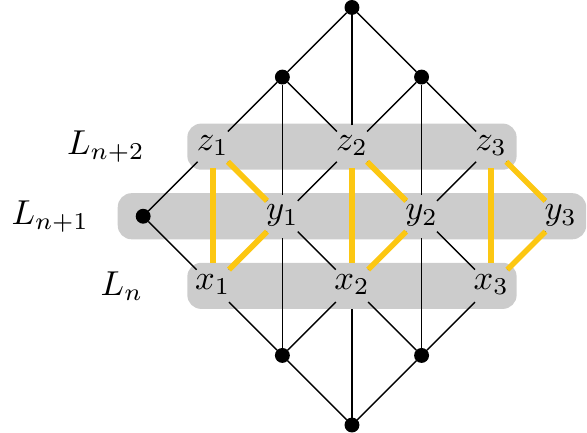}
    \caption{
      Three triangles in~$\Grid_4 \subseteq N_4$ with vertices in levels~$L_4, L_5, $ and $ L_6$.
    }
    \label{fig:triangles-in-grid}
  \end{figure}

  We now define an ordering $<_T$ on the triangles and use it to find a $(p + 1)$\nobreakdash-twist.
  We define~$T_i <_T T_j$ if and only if~$x_i < x_j$.
  A subsequence $y_{i_1}, \ldots, y_{i_s}$ of $y_1,\ldots,y_r$ is \emph{increasing} if its ordering corresponds to~$<_T$, that is if~$T_{i_1} <_T \cdots <_T T_{i_s}$.
  Similarly, a subsequence $y_{i_1}, \ldots, y_{i_s}$ of $y_1,\ldots,y_r$ is called \emph{decreasing} if their reverse ordering corresponds to~$<_T$, that is if~$T_{i_1} >_T \cdots >_T T_{i_s}$.
  Increasing and decreasing subsequences of~$z_1, \ldots, z_r$ in level~$L_{n + 2}$ are defined analogously.

  We now only consider the subgraph of~$N_n$ that is given by the triangles~$T_1, \ldots, T_r$.
  That is, a neighbor of a vertex~$v$ refers to a vertex in the same triangle as~$v$.
  If there is an increasing subsequence of $y_1,\ldots,y_r$ or of $z_1,\ldots,z_r$ of length~$p + 1$, then we have a~$(p + 1)$\nobreakdash-twist between these vertices and their neighbors in~$L_n$.
  Hence, the longest increasing subsequences of $y_1,\ldots,y_r$ and $z_1,\ldots,z_r$ have length at most~$p$.
  By the Erd\H{o}s-Szekeres theorem, there exists a decreasing subsequence $y_{i_1},\ldots,y_{i_s}$ of $y_1,\ldots,y_r$ of length~$s = p^2 + 1$.
  Again by the Erd\H{o}s-Szekeres theorem, there exists a decreasing subsequence $z_{i'_1},\ldots,z_{i'_t}$ of $z_{i_1},\ldots,z_{i_s}$ of length~$t = p + 1$.
  But then $y_{i'_t} < \cdots < y_{i'_1} < z_{i'_t} < \cdots < z_{i'_1}$ form a~$(p + 1)$\nobreakdash-twist as $L_{n + 1} < L_{n + 2}$.
\end{proof}

We conclude by \cref{lem:G5-ordering} that for~$p = 4$ and~$n = p^3 + 2 = 66$ every level-separating topological ordering of $N_n$ contains a~$5$\nobreakdash-twist.
Further, by \cref{lem:separating-levels-N}, there is an~$n' \geq n$ such that every topological ordering~$<$ of~$N_{n',5}^\ast$ contains a copy of~$N_n$ whose levels are separated by~$<$ (i.e., $n' = n + 2(n-1) = 192$ as in the proof).
Together this yields~$\pn(N_{n'}) \geq \tn(N_{n'}) \geq 5$, proving \cref{main:5-twist}.

Finally, we remark that N-grids have bounded page number but it is not obvious whether five pages suffice for all N-grids.
However, separating the levels of N-grids works only with 5-fences, which is why new ideas are needed for any significant improvement.

\section{Conclusions}%
\label{sec:xonclusions}

In this paper, we improve both the lower and the upper bound on the maximum page number among upward planar graphs.
Concerning the lower bound, we remark that \cref{lem:G5-ordering} does not depend on the size of the twist to be enforced but yields arbitrarily large twists. 
That is, for pushing the lower bound further it suffices to find a large enough upward planar graph whose vertices can be partitioned into levels, i.e., into sets of vertices that are separated by any topological ordering with no large twist.
However, it is crucial that there are edges connecting non-consecutive levels.
We also expect the concept of fences to prove useful for improving the lower bound further as we only need to consider topological orderings that respect the augmented edges.

The main contribution of this paper is the first sublinear upper bound on the page number of upward planar graphs in terms of their number of vertices.
We remark that when applying \cref{lem:width} repeatedly, many edges are embedded multiple times. 
In fact, we only need the edges of $ G'[X] $ in the last application of the lemma, whereas we use the embedding of the edges in $ E_\Delta $ in all rounds.
In light of this observation, we see potential improvements in reducing the number of pages needed for~$ E_\Delta $ (at the expense of the number of pages needed for $ G'[X] $) or in reducing the number of applications of \cref{lem:width} (e.g., by covering the edges of $ E_\Delta $ with \cref{lem:height}).
Both would lead to an upper bound of $ \mathcal{O}(\sqrt{n \log(n)}) $.
To improve the bound beyond that, we think that new approaches are necessary.


In addition to the sublinear upper bound, we attack the problem of bounding the page number of upward planar graphs by showing that families of upward planar graphs with bounded width or bounded height have bounded page number.
However, the initial question by Nowakowski and Parker~\cite{NP89} whether planar posets, and more generally upward planar graphs, have bounded page number, still remains open.

\bibliography{references}

\end{document}